\documentclass[a4paper,12pt,final]{amsart}
\usepackage{times,a4wide,mathrsfs,amssymb,amsmath,amsthm,enumerate,xypic,tikzsymbols,dsfont}

\newcommand{\C}{\mathbb{C}}

\newcommand{\QQ}{\mathbb{Q}}
\newcommand{\NN}{\mathbb{N}}
\newcommand{\PP}{\mathbb{P}}
\newcommand{\SSS}{\mathbb{S}}

\newcommand{\OO}{\mathcal O}
\newcommand{\Ss}{\mathcal S}

\newcommand{\XX}{\mathcal X}

\newcommand{\HH}{\mathcal H}
\newcommand{\NNN}{\mathcal N}

\newcommand{\MM}{\mathcal M}

\newcommand{\wt}{\widetilde}

\newcommand{\one}{\mathds{1}}

\DeclareMathOperator{\ima}{Im}

\DeclareMathOperator{\rank}{rank}

\DeclareMathOperator{\Gr}{Gr}
\DeclareMathOperator{\SGr}{SGr}
\DeclareMathOperator{\OGr}{OGr}
\DeclareMathOperator{\Dec}{Dec}

\newtheorem{theorem}{Theorem}[section]

\newtheorem{corollary}[theorem]{Corollary}
\newtheorem{proposition}[theorem]{Proposition}
\newtheorem{conjecture}[theorem]{Conjecture}
\newtheorem{remark}[theorem]{Remark}
\newtheorem{definition}[theorem]{Definition}
\newtheorem{example}[theorem]{Example}
\newtheorem{convention}{Conventions}

\newtheorem{notation}[theorem]{Notation}

\newtheorem{nonumbering}{Theorem}

\newtheorem{nonumberingc}{Corollary}

\newtheorem{nonumberingt}{Acknowledgements}

\begin{document}

\author[Robert Laterveer]
{Robert Laterveer}

\address{Institut de Recherche Math\'ematique Avanc\'ee,
CNRS -- Universit\'e 
de Strasbourg,\
7 Rue Ren\'e Des\-car\-tes, 67084 Strasbourg CEDEX,
FRANCE.}
\email{laterv@math.unistra.fr}

\title{On the Chow ring of Fano varieties on the Fatighenti-Mongardi list}

\begin{abstract} Conjecturally, Fano varieties of K3 type admit a multiplicative Chow--K\"unneth decomposition, in the sense of Shen--Vial. We prove this for many of the families of Fano varieties of K3 type constructed by Fatighenti--Mongardi.
This has interesting consequences for the Chow ring of these varieties.
%and also for a certain tautological subring of the Chow ring of powers of these Fano varieties.
 \end{abstract}

\thanks{Supported by ANR grant ANR-20-CE40-0023.}

%\keywords{Algebraic cycles, Chow groups, motive, Bloch--Beilinson filtration, Beauville's ``splitting property'' conjecture, multiplicative Chow--K\"unneth decomposition, Fano variety, K3 surface}
%\subjclass[2010]{Primary 14C15, 14C25, 14C30.}

\maketitle

\section{Introduction}

Given a smooth projective variety $Y$ over $\C$, let $A^i(Y):=CH^i(Y)_{\QQ}$ denote the Chow groups of $Y$ (i.e. the groups of codimension $i$ algebraic cycles on $Y$ with $\QQ$-coefficients, modulo rational equivalence). The intersection product defines a ring structure on $A^\ast(Y)=\bigoplus_i A^i(Y)$, the Chow ring of $Y$ \cite{F}. In the case of K3 surfaces, this ring structure has a remarkable property:

\begin{theorem}[Beauville--Voisin \cite{BV}]\label{K3} Let $S$ be a K3 surface. 
%Let $D_i, D_i^\prime\in A^1(S)$ be a finite number of divisors. Then
%  \[ \sum_i D_i\cdot D_i^\prime=0\ \ \ \hbox{in}\ A^2(S)_{}\ \Leftrightarrow\ \sum_i D_i\cdot D_i^\prime=0\ \ \ \hbox{in}\ H^4(S,\QQ)\ .\]
The $\QQ$-subalgebra
  \[  R^\ast(S):=  \bigl\langle  A^1(S), c_j(S) \bigr\rangle\ \ \ \subset\ A^\ast(S) \]
  injects into cohomology under the cycle class map.
  \end{theorem}

Motivated by the cases of K3 surfaces and abelian varieties, Beauville \cite{Beau3} has conjectured that for certain special varieties, the Chow ring should admit a multiplicative splitting. To make concrete sense of Beauville's elusive ``splitting property conjecture'', Shen--Vial \cite{SV} have introduced the concept of {\em multiplicative Chow--K\"unneth decomposition\/}. It seems both interesting and difficult to better understand the class of special varieties admitting such a decomposition.

In \cite{S2}, the following conjecture is raised:

\begin{conjecture}\label{conj} Let $X$ be a smooth projective Fano variety of K3 type (i.e. $\dim X=2d$ and the Hodge numbers $h^{p,q}(X)$ are $0$ for all $p\not=q$ except for $h^{d-1,d+1}(X)=h^{d+1,d-1}(X)=1$). Then $X$ has a multiplicative Chow--K\"unneth decomposition.
\end{conjecture}

This conjecture is verified in some special cases \cite{37}, \cite{39}, \cite{40}, \cite{FLV2}, \cite{S2}.
This paper aims to contribute to this program. The main result is as follows:

\begin{nonumbering}[=Theorem \ref{main}] Let $X$ be a smooth Fano variety in one of the families of Table \ref{table:1}. Then $X$ has a multiplicative Chow--K\"unneth decomposition. 
%The Chern classes $c_j(X)$ lie in $A^j_{(0)}(X)$. The correspondence $\Gamma$ linking $X$ and the associated K3 surface lies in $A^\ast_{(0)}(X\times S)$.
\end{nonumbering}

Table \ref{table:1} lists Fano varieties $X$ of K3 type that were constructed by Fatighenti--Mongardi \cite{FM} as hypersurfaces in products of Grassmannians. The K3 surfaces $S$ in Table \ref{table:1} are shown in \cite{FM} to be associated to $X$ on the level of Hodge theory, and on the level of derived categories. In some cases, the geometric relation between $X$ and $S$ is straightforward (e.g., for B1 and B2 the Fano variety $X$ is a blow-up with center the K3 surface $S$); in other cases the geometric relation is more indirect (e.g. for M1, M6, M7, M8, M9, M10 the Fano variety $X$ is related to the K3 surface $S$ via the so-called ``Cayley's trick'', cf. \cite{FM} and subsection \ref{ss:cay} below).

To prove Theorem \ref{main}, we have devised a general criterion (Proposition \ref{crit}), which we hope might apply to other Fano varieties of K3 type. To verify the criterion, one needs a motivic relation between the Fano variety $X$ and the associated K3 surface $S$, and one needs a certain instance of the {\em Franchetta property\/}.

\bigskip
\bigskip
\begin{table}[h]
\centering
\begin{tabular}{||c c c c c c||} 
 \hline
 $\stackrel{\hbox{Label}}{\hbox{in\ \cite{FM}}}$ & $X\subset U$ & $\dim X$ & $\rho(X)$   & $\stackrel{\hbox{Genus of}}{\hbox{associated\ K3}}$  &$\stackrel{\hbox{Also}}
 {\hbox{occurs\ in}}$ \\ 
 [0.5ex] 
 \hline\hline
 B1 & $X_{(2,1,1)}\subset\PP^3\times\PP^1\times\PP^1$ & 4 & 3  & 7& \cite{40}\\ 
  B2 & $X_{(2,1)}\subset\Gr(2,4)\times\PP^1$ & 4 & 3  & 5& \cite{40}\\ 
   M1 & $X_{(1,1,1)}\subset\PP^3\times\PP^3\times\PP^3$ & 8 & 3  & 3& \cite{IM}\\ 
 M3 & $X_{(1,1)}\subset\Gr(2,5)\times Q_5$ & 10 & 2  & 6&\\
 M4 & $X_{(1,1)}\subset\SGr(2,5)\times Q_4$  & 8 & 2 & 6 &\\
 M6 & $X_{(1,1)}\subset \SSS_5\times\PP^7$ & 16 & 2 & 7&\\
 M7 &   $X_{(1,1)}\subset\Gr(2,6)\times \PP^5$ & 12 & 2 & 8 &\\ 
 M8 &  $X_{(1,1)}\subset\SGr(2,6)\times \PP^4$ & 10 & 2 & 8 &\\ 
 M9 & $X_{(1,1)}\subset S_2 \Gr(2,6)\times \PP^3$ & 8 & 2  & 8&\\  
 M10 & $X_{(1,1)}\subset\SGr(3,6)\times \PP^3$ & 8 & 2 & 9&\\  
% M13 &  $X_{(1,1)}\subset\Gr(2,8)\times \PP^3$ & 14 & 2  & 3&\\  
% S1 & $X_{(1,1,1,1)}\subset \Gr(2,8)$ & 8 & 1 & 3& \cite{ST}\\
 S2 &  $X_{1}\subset \OGr(2,8)$ & 8 & 2 & 7& \cite{S2}\\ 
 [1ex] 
 \hline
\end{tabular}
\caption{Families of Fano varieties of K3 type. (As in \cite{FM}, $\Gr(k,m)$ denotes the Grassmannian of $k$-dimensional subspaces of an $m$-dimensional vector space. $\SGr(k,m)$, $S_2 \Gr(k,m)$ and $\OGr(k,m)$ denote the symplectic resp. bisymplectic resp. orthogonal Grassmannian. $\SSS_5$ denotes a connected component of $\OGr(5,10)$, and $Q_m$ is an $m$-dimensional smooth quadric.)} 
%The label refers to the notation of \cite{FM}.}
\label{table:1}
\end{table}
\medskip

As a consequence of our main result, the Chow ring of these Fano varieties behaves like the Chow ring of a K3 surface:

\begin{nonumberingc}[=Corollary \ref{cor}] Let $X\subset U$ be the inclusion of a Fano variety $X$ in its ambient space $U$, where $X,U$ are as in Table \ref{table:1}. Let $\dim X=2d$. Let $R^\ast(X)\subset A^\ast(X)$ be the $\QQ$-subalgebra
  \[ R^\ast(X):=\Bigl\langle A^1(X), A^2(X), \ldots, A^d(X), c_j(X),\ima\bigl(A^\ast(U)\to A^\ast(X)\bigr)\Bigr\rangle\ \ \ \subset A^\ast(X)\ .\]
  Then $R^\ast(X)$ injects into cohomology under the cycle class map.
 \end{nonumberingc}

We end this introduction with a challenge. Fatighenti--Mongardi have constructed some more Fano varieties of K3 type for which it would be nice to settle Conjecture \ref{conj} (in particular the families labelled M13 and S1 in \cite{FM}, for which I have not been able to check condition (c3) or (c3$^\prime$) of the general criterion Proposition \ref{crit}).

Additionally, the following are some Fano varieties of K3 type in the litterature for which Conjecture \ref{conj} is still open, and for which 
the methods of the present paper are not sufficiently strong:
%we have not been able to verify our general criterion (Proposition \ref{crit}): 
K\"uchle fourfolds of type $c5$, Pl\"ucker hyperplane sections of $\Gr(3,10)$, intersections of $\Gr(2,8)$ with 4 Pl\"ucker hyperplanes,
Gushel--Mukai fourfolds and sixfolds. It would be interesting to 
%find means to verify our criterion for these families.
devise new methods to treat these families.

 \vskip0.6cm

\begin{convention} In this article, the word {\sl variety\/} will refer to a reduced irreducible scheme of finite type over $\C$. A {\sl subvariety\/} is a (possibly reducible) reduced subscheme which is equidimensional. 

{\bf All Chow groups will be with rational coefficients}: we denote by $A_j(Y)$ the Chow group of $j$-dimensional cycles on $Y$ with $\QQ$-coefficients; for $Y$ smooth of dimension $n$ the notations $A_j(Y)$ and $A^{n-j}(Y)$ are used interchangeably. 
The notation $A^j_{hom}(Y)$ will be used to indicate the subgroup of homologically trivial cycles.
%For a morphism $f\colon X\to Y$, we will write $\Gamma_f\in A_\ast(X\times Y)$ for the graph of $f$.

The contravariant category of Chow motives (i.e., pure motives with respect to rational equivalence as in \cite{Sc}, \cite{MNP}) will be denoted 
$\MM_{\rm rat}$.
%We will write $H^j(X)$ to indicate singular cohomology $H^j(X,\QQ)$.
\end{convention}

\section{Preliminaries}

\subsection{MCK decomposition}
\label{ss:mck}

\begin{definition}[Murre \cite{Mur}] Let $X$ be a smooth projective variety of dimension $n$. We say that $X$ has a {\em CK decomposition\/} if there exists a decomposition of the diagonal
   \[ \Delta_X= \pi^0_X+ \pi^1_X+\cdots +\pi_X^{2n}\ \ \ \hbox{in}\ A^n(X\times X)\ ,\]
  such that the $\pi^i_X$ are mutually orthogonal idempotents and $(\pi_X^i)_\ast H^\ast(X,\QQ)= H^i(X,\QQ)$.
  
  (NB: ``CK decomposition'' is shorthand for ``Chow--K\"unneth decomposition''.)
\end{definition}

\begin{remark} The existence of a CK decomposition for any smooth projective variety is part of Murre's conjectures \cite{Mur}, \cite{J4}. 
\end{remark}

\begin{definition}[Shen--Vial \cite{SV}] Let $X$ be a smooth projective variety of dimension $n$. Let $\Delta_X^{sm}\in A^{2n}(X\times X\times X)$ be the class of the small diagonal
  \[ \Delta_X^{sm}:=\bigl\{ (x,x,x)\ \vert\ x\in X\bigr\}\ \subset\ X\times X\times X\ .\]
  An {\em MCK decomposition\/} is a CK decomposition $\{\pi_X^i\}$ of $X$ that is {\em multiplicative\/}, i.e. it satisfies
  \[ \pi_X^k\circ \Delta_X^{sm}\circ (\pi_X^i\times \pi_X^j)=0\ \ \ \hbox{in}\ A^{2n}(X\times X\times X)\ \ \ \hbox{for\ all\ }i+j\not=k\ .\]
  
 (NB: ``MCK decomposition'' is shorthand for ``multiplicative Chow--K\"unneth decomposition''.) 
  
% A {\em weak MCK decomposition\/} is a CK decomposition $\{\pi^X_i\}$ of $X$ that satisfies
%    \[ \Bigl(\pi^X_k\circ \Delta_X^{sm}\circ (\pi^X_i\times \pi^X_j)\Bigr){}_\ast (a\times b)=0 \ \ \ \hbox{for\ all\ } a,b\in\ A^\ast(X)\ .\]
  \end{definition}
  
  \begin{remark} The small diagonal (seen as a correspondence from $X\times X$ to $X$) induces the {\em multiplication morphism\/}
    \[ \Delta_X^{sm}\colon\ \  h(X)\otimes h(X)\ \to\ h(X)\ \ \ \hbox{in}\ \MM_{\rm rat}\ .\]
 Let us assume $X$ has a CK decomposition
  \[ h(X)=\bigoplus_{i=0}^{2n} h^i(X)\ \ \ \hbox{in}\ \MM_{\rm rat}\ .\]
  By definition, this decomposition is multiplicative if for any $i,j$ the composition
  \[ h^i(X)\otimes h^j(X)\ \to\ h(X)\otimes h(X)\ \xrightarrow{\Delta_X^{sm}}\ h(X)\ \ \ \hbox{in}\ \MM_{\rm rat}\]
  factors through $h^{i+j}(X)$.
  
  If $X$ has an MCK decomposition, then setting
    \[ A^i_{(j)}(X):= (\pi_X^{2i-j})_\ast A^i(X) \ ,\]
    one obtains a bigraded ring structure on the Chow ring: that is, the intersection product sends $A^i_{(j)}(X)\otimes A^{i^\prime}_{(j^\prime)}(X) $ to  $A^{i+i^\prime}_{(j+j^\prime)}(X)$.
    
      It is expected that for any $X$ with an MCK decomposition, one has
    \[ A^i_{(j)}(X)\stackrel{??}{=}0\ \ \ \hbox{for}\ j<0\ ,\ \ \ A^i_{(0)}(X)\cap A^i_{hom}(X)\stackrel{??}{=}0\ ;\]
    this is related to Murre's conjectures B and D, that have been formulated for any CK decomposition \cite{Mur}.

  The property of having an MCK decomposition is restrictive, and is closely related to Beauville's ``splitting property' conjecture'' \cite{Beau3}. 
  To give an idea: hyperelliptic curves have an MCK decomposition \cite[Example 8.16]{SV}, but the very general curve of genus $\ge 3$ does not have an MCK decomposition \cite[Example 2.3]{FLV2}. As for surfaces: a smooth quartic in $\PP^3$ has an MCK decomposition, but a very general surface of degree $ \ge 7$ in $\PP^3$ should not have an MCK decomposition \cite[Proposition 3.4]{FLV2}.
For more detailed discussion, and examples of varieties with an MCK decomposition, we refer to \cite[Section 8]{SV}, as well as \cite{V6}, \cite{SV2}, \cite{FTV}, \cite{37}, \cite{39}, \cite{40}, \cite{S2}, \cite{46}, \cite{38}, \cite{FLV2}.
   \end{remark}

 \subsection{The Franchetta property}
 \label{ss:fr}

 \begin{definition} Let $\XX\to B$ be a smooth projective morphism, where $\XX, B$ are smooth quasi-projective varieties. We say that $\XX\to B$ has the {\em Franchetta property in codimension $j$\/} if the following holds: for every $\Gamma\in A^j(\XX)$ such that the restriction $\Gamma\vert_{X_b}$ is homologically trivial for the very general $b\in B$, the restriction $\Gamma\vert_b$ is zero in $A^j(X_b)$ for all $b\in B$.
 
 We say that $\XX\to B$ has the {\em Franchetta property\/} if $\XX\to B$ has the Franchetta property in codimension $j$ for all $j$.
 \end{definition}
 
 This property is studied in \cite{PSY}, \cite{BL}, \cite{FLV}, \cite{FLV3}.
 
 \begin{definition} Given a family $\XX\to B$ as above, with $X:=X_b$ a fiber, we write
   \[ GDA^j_B(X):=\ima\Bigl( A^j(\XX)\to A^j(X)\Bigr) \]
   for the subgroup of {\em generically defined cycles}. 
  In a context where it is clear to which family we are referring, the index $B$ will often be suppressed from the notation.
  \end{definition}
  
  With this notation, the Franchetta property amounts to saying that $GDA^\ast_B(X)$ injects into cohomology, under the cycle class map.

\subsection{Cayley's trick and motives}
\label{ss:cay}

\begin{theorem}[Jiang \cite{Ji}]\label{ji} Let $ E\to U$ be a vector bundle of rank $r\ge 2$ over a smooth projective variety $U$, and let $S:=s^{-1}(0)\subset U$ be the zero locus of a regular section $s\in H^0(U,E)$ such that $S$ is smooth of dimension $\dim U-\rank E$. Let $X:=w^{-1}(0)\subset \PP(E)$ be the zero locus of the regular section $w\in H^0(\PP(E),\OO_{\PP(E)}(1))$ that corresponds to $s$ under the natural isomorphism $H^0(U,E)\cong H^0(\PP(E),\OO_{\PP(E)}(1))$, and assume $X$ is smooth. There is an isomorphism of Chow motives
    \[   h(X)\cong h(S)(1-r)\oplus \bigoplus_{i=0}^{r-2} h(U)(-i)\ \ \ \hbox{in}\ \MM_{\rm rat}\ .\]
    \end{theorem}
    
    \begin{proof} This is \cite[Corollary 3.2]{Ji}, which more precisely gives an isomorphism of {\em integral\/} Chow motives. For later use, we now give some details about the isomorphism as constructed in loc. cit.. Let 
    \[ \Gamma:= X\times_U S\ \ \subset\ X\times S \]
  (this is equal to $\PP(\NNN_i)=\HH_s\times_X Z$ in the notation of loc. cit.).  Let
    \[ \Pi_i\ \ \in\ A^\ast(X\times U)\ \ \ \ (i=0, \ldots, r-2) \]
    be correspondences inducing the maps $ (\pi_i)_\ast$ of loc. cit., i.e.  
       \[ (\Pi_i)_\ast= (\pi_i)_\ast :=  (q_{i+1})_\ast \iota_\ast\colon\ \ A^j(X)\ \to\ A^{j-i}(U)\ ,\]
       where $\iota\colon X\hookrightarrow\PP(E)$ is the inclusion morphism, and the $(q_{i+1})_\ast\colon A_\ast(\PP(E))\to A_\ast(U)$ are defined in loc. cit. in terms of the projective bundle formula for $q\colon E\to U$. As indicated in \cite[Corollary 3.2]{Ji} (cf. also \cite[text preceeding Corollary 3.2]{Ji}), there is an isomorphism
       \[   \Bigl( \Gamma, \Pi_0,\Pi_1,\ldots, \Pi_{r-2}\Bigr)\colon\ \ h(X)\ \xrightarrow{\cong}\ h(S)(1-r)\oplus \bigoplus_{i=0}^{r-2} h(U)(-i)\ \ \ \hbox{in}\ \MM_{\rm rat}\ .\]    
%       The inverse isomorphism is given as
%       \[  
    \end{proof}

\begin{remark} In the set-up of Theorem \ref{ji}, a cohomological relation between $X$ and $S$ was established in \cite[Prop. 4.3]{Ko} (cf. also \cite[section 3.7]{IM0}, as well as \cite[Proposition 46]{BFM} for a generalization). A relation on the level of derived categories was established in \cite[Theorem 2.10]{Or} (cf. also \cite[Theorem 2.4]{KKLL} and \cite[Proposition 47]{BFM}).
\end{remark}

 We now make the natural observation that the isomorphism of Theorem \ref{ji} behaves well with respect to families, in the following sense: 
 
 \begin{notation} Let $X, S, U$ and $E\to U$ be as in Theorem \ref{ji}. Let $B\subset\PP H^0(\PP(E),\OO_{\PP(E)}(1))$ be the Zariski open such that 
 both $X:=X_b\subset\PP(E)$ and $S:=S_b\subset U$ are smooth of the expected dimension. Let
    \[  \XX\to B\ ,\ \ \ \Ss\to B \]
    denote the universal families. 
%    
%    Given $m\in\NN$, we will write
%    \[ \XX^{m/B}:=  \XX\times_B \XX\times_B \cdots \times_B \XX \]
%    for the $m$-fold fiber product.
     \end{notation}

 \begin{proposition}\label{ji2} Let $X, S, U$ be as in Theorem \ref{ji}. Assume $U$ has trivial Chow groups. For any $m\in\NN$, there are injections
    \[  GDA^j(X^m)\ \hookrightarrow\ GDA^{j+m-mr}(S^m)\oplus \bigoplus GDA^\ast(S^{m-1})\oplus \cdots \cdots \oplus \QQ^s\ .\]
       \end{proposition} 
       
       \begin{proof} (NB: we will not really need this proposition below, but we include it because it makes some arguments easier, cf. footnote 1 below.)
       
       Let us first do the case $m=1$. The isomorphism of Theorem \ref{ji} is {\em generically defined\/}, i.e. there exist relative correspondences
       $\Gamma_B,\Pi_i^B$ fitting into a commutative diagram
       \begin{equation}\label{dia}  \begin{array}[c]{ccc}     A^j(\XX) & \xrightarrow{\bigl( (\Gamma_B)_\ast,(\Pi_0^B)_\ast,\ldots,(\Pi^B_{r-2})_\ast\bigr)} & A^{j+1-r}(\Ss)\oplus \bigoplus_{i=0}^{r-2} A^{j-i}(U\times B)\\
       &&\\
       \downarrow&&\downarrow\\
       &&\\
        A^j(X) & \xrightarrow{\bigl( \Gamma_\ast,(\Pi_0)_\ast,\ldots,(\Pi_{r-2})_\ast\bigr)} & \ A^{j+1-r}(S)\oplus \bigoplus_{i=0}^{r-2} A^{j-i}(U),\\ 
        \end{array}\end{equation}
        where vertical arrows are restrictions to a fiber, and the lower horizontal arrow is the isomorphism of Theorem \ref{ji}.
        Indeed, $\Gamma_B$ can be defined as
          \[ \Gamma_B:= \XX\times_{U\times B}\Ss\ \ \subset\ \XX\times_B \Ss\ .\]
          The $\Pi_i$ are also generically defined (just because the graph of the embedding $\iota\colon X\hookrightarrow \PP(E)$ is generically defined). This gives
          relative correspondences $\Gamma_B, \Pi_i^B$ over $B$ such that the restriction to a fiber over $b\in B$ gives back the correspondences $\Gamma,\Pi_i$ of Theorem \ref{ji}. The fact that this makes diagram \eqref{dia} commute is \cite[Lemma 8.1.6]{MNP}.
          The commutative diagram \eqref{dia} implies that there is an injective map
          \begin{equation}\label{1} GDA^j(X)\ \hookrightarrow\ GDA^{j+1-r}(S)\oplus \bigoplus A^\ast(U) = GDA^{j+1-r}(S)\oplus \QQ^s\ .\end{equation}
%   For the surjectivity (which we will not use below), we note that Voisin's Hilbert schemes argument  \cite[Proposition 3.7]{V0} (in the precise form used in \cite[Proposition 2.11]{Lacub}) implies that the inverse to the isomorphism of Theorem \ref{ji} is also generically defined (alternatively, one may check directly that the inverse isomorphism in \cite[Corollary 3.2]{Ji} is generically defined).
%   Using a diagram similar to \eqref{dia} (with horizontal arrows in the opposite direction), this gives that the map \eqref{1} is also surjective and so the $m=1$ case is proven.
   
 The argument for $m>1$ is similar: the isomorphism of motives of Theorem \ref{ji}, combined with the fact that $U$ has trivial Chow groups (and so $h(U)\cong \oplus \one(\ast)$) induces an isomorphism of Chow groups
   \begin{equation}\label{2}  A^j(X^m)\   \xrightarrow{\cong}\ A^{j+m-mr}(S^m)\oplus \bigoplus A^\ast(S^{m-1})\oplus \cdots \cdots \oplus \QQ^s\ .\end{equation}
   Here the map from left to right is given by various combinations of the correspondences $\Gamma$ and $\Pi_i$. As we have seen these correspondences are generically defined, and so their products are also generically defined. It follows as above that the map \eqref{2} preserves generically defined cycles.         
         \end{proof}

  \subsection{A Franchetta-type result}
  
  \begin{proposition}\label{spread} Let $Y$ be a smooth projective variety with trivial Chow groups (i.e. $A^\ast_{hom}(Y)=0$). Let $L_1,\ldots,L_r\to Y$ be very ample line bundles, and let
  $\XX\to B$ be the universal family of smooth complete intersections of type $X=Y\cap H_1\cap\cdots\cap H_r$, where $H_j\in\vert L_j\vert$.
  Assume the fibers $X$ have $H^{\dim X}_{tr}(X,\QQ)\not=0$.
  There is inclusion
    \[ \ker \Bigl( GDA^{\dim X}_B(X\times X)\to H^{2\dim X}(X\times X,\QQ)\Bigr)\ \ \subset\ \Bigl\langle (p_1)^\ast GDA^\ast_B(X), (p_2)^\ast GDA^\ast_B(X)  \Bigr\rangle\ .\]
   \end{proposition}
   
   \begin{proof} This is essentially equivalent to Voisin's ``spread'' result \cite[Proposition 1.6]{V1} (cf. also \cite[Proposition 5.1]{LNP} for a reformulation). For completeness, we include a quick proof. Let $\bar{B}:=\PP H^0(Y,L_1\oplus\cdots\oplus L_r)$ (so that $B\subset\bar{B}$ is a Zariski open), and let 
   us consider the projection
   \[ \pi\colon\ \  \XX\times_{\bar{B}} \XX\ \to\ Y\times Y\ .\]
   Using the very ampleness assumption, one finds that $\pi$ is a $\PP^s$-bundle over $(Y\times Y)\setminus \Delta_Y$, and a $\PP^t$-bundle over $\Delta_Y$.
   That is, $\pi$ is what is termed a {\em stratified projective bundle\/} in \cite{FLV}. As such, \cite[Proposition 5.2]{FLV} implies the equality
      \begin{equation}\label{stra}  GDA^\ast_B(X\times X)= \ima\Bigl( A^\ast(Y\times Y)\to A^\ast(X\times X)\Bigr) +  \Delta_\ast GDA_B^\ast(X)\ ,\end{equation}
      where $\Delta\colon X\to X\times X$ is the inclusion along the diagonal. Since $Y$ has trivial Chow groups, one has $A^\ast(Y\times Y)\cong A^\ast(Y)\otimes A^\ast(Y)$.
       Base-point freeness of the $L_j$ implies $\XX\to Y$ has the structure of a projective bundle; it is then readily seen (by a direct argument or by simply applying once more \cite[Proposition 5.2]{FLV}) that 
        \[  GDA^\ast_B(X)=\ima\bigl( A^\ast(Y)\to A^\ast(X)\bigr)\ .\]
       The equality \eqref{stra} thus reduces to      
      \[ GDA^\ast_B(X\times X)=\Bigl\langle (p_1)^\ast GDA^\ast_B(X), (p_2)^\ast GDA^\ast_B(X), \Delta_X\Bigr\rangle\ \]   
      (where $p_1, p_2$ denote the projection from $X\times X$ to first resp. second factor). The assumption that $X$ has non-zero transcendental cohomology
      implies that the class of $\Delta_X$ is not decomposable in cohomology. It follows that
      \[ \begin{split}  \ima \Bigl( GDA^{\dim X}_B(X\times X)\to H^{2\dim X}(X\times X,\QQ)\Bigr) =&\\
       \ima\Bigl(  \Dec^{\dim X}(X\times X)\to H^{2\dim X}(X\times X,\QQ)\Bigr)& \oplus \QQ[\Delta_X]\ ,\\
       \end{split}\]
      where we use the shorthand 
       \[ \Dec^j(X\times X):= \Bigl\langle (p_1)^\ast GDA^\ast_B(X), (p_2)^\ast GDA^\ast_B(X)\Bigr\rangle\cap A^j(X\times X) \ \]     
       for the {\em decomposable cycles\/}. 
       We now see that if $\Gamma\in GDA^{\dim X}(X\times X)$ is homologically trivial, then $\Gamma$ does not involve the diagonal and so $\Gamma\in \Dec^{\dim X}(X\times X)$.
       This proves the proposition.
         \end{proof}
  
  \begin{remark} Proposition \ref{spread} has the following consequence: if the family $\XX\to B$ has the Franchetta property, then $\XX\times_B \XX\to B$ has the Franchetta property in codimension $\dim X$.
   \end{remark}

 \subsection{HPD and motives}
 
 \begin{theorem}\label{hpd} Let $Y_1, Y_2\subset\PP(V)$ be smooth projective varieties with trivial Chow groups (i.e. $A^\ast_{hom}(Y_j)=0$), and let $Y_2^\vee\subset\PP(V^\vee)$ be the HPD dual of $Y_2$. Let $H\subset \PP(V)\times\PP(V)$ be a $(1,1)$-divisor, and let $f_H\colon \PP(V)\to\PP(V^\vee)$ be the morphism defined by $H$. Assume that the varieties
   \[ \begin{split}  X&:= (Y_1\times Y_2)\cap H\ ,\\
                                  S&:=Y_1\cap (f_H)^{-1}(Y_2^\vee)\\
                                  \end{split}\]
    are smooth and dimensionally transverse. Assume moreover that the Hodge conjecture holds for $S$, that $H^j(S,\QQ)$ is algebraic for $j\not=\dim S$ and that $H^{\dim S}(S,\QQ)$ is not completely algebraic.
    Then there is a split injection of Chow motives
    \[ h(X)\ \hookrightarrow\ h(S)(-m)\oplus \bigoplus\one(\ast)\ \ \ \hbox{in}\ \MM_{\rm rat}    \ ,\]
    where $m:={1\over 2}(\dim X-\dim S)$.        
    
    (In particular, one has vanishing
       \[    A^j_{hom}(X) =0\ \ \ \forall\ j  > {1\over 2}(\dim X+\dim S)\ .)\]                  
      \end{theorem} 
  
  \begin{proof} Using the HPD formalism, it is proven in \cite[Proposition 2.4]{FM} that there exists a semi-orthogonal decomposition
    \begin{equation}\label{so} D^b(X)=\bigl\langle D^b(S), A_1,\ldots, A_s\bigr\rangle\ ,\end{equation}
    where the $A_j$ are some exceptional objects. Using Hochschild homology and the Kostant--Rosenberg isomorphism (cf. for instance \cite[Sections 1.7 and 2.5]{Kuz}), this implies that there exist correspondences $\Phi^\prime$ and $\Xi^\prime$ such that
      \[  H^{\ast}_{tr}(X,\QQ)\ \xrightarrow{(\Phi^\prime)_\ast}\ H^{\ast}_{tr}(S,\QQ)\ \xrightarrow{(\Xi^\prime)_\ast}\   H^{\ast}_{tr}(X,\QQ) \]
    is the identity. (Here, $H^\ast_{tr}( ,\QQ)$ denotes the orthogonal complement of the algebraic part of cohomology.)
    By assumption $H^\ast_{tr}(S,\QQ)=H^{\dim S}_{tr}(S,\QQ)$, and by weak Lefschetz $H^\ast_{tr}(X,\QQ)=  H^{\dim X}_{tr}(X,\QQ)$, and so we actually have that 
      \[  H^{\dim X}_{tr}(X,\QQ)\ \xrightarrow{(\Phi^\prime)_\ast}\ H^{\dim S}_{tr}(S,\QQ)\ \xrightarrow{(\Xi^\prime)_\ast}\   H^{\dim X}_{tr}(X,\QQ) \]
    is the identity. Again using Hochschild homology and the Kostant--Rosenberg isomorphism, we see that the Hodge conjecture for $S$, plus the decomposition \eqref{so}, implies the Hodge conjecture for $X$. This means that we can find correspondences $\Phi$ and $\Xi$ such that
       \[  H^{\ast}_{}(X,\QQ)\ \xrightarrow{\Phi_\ast}\ H^{\dim S}_{}(S,\QQ) \oplus \bigoplus \QQ(-j)\ \xrightarrow{\Xi_\ast}\   H^{\ast}_{}(X,\QQ) \]
    is the identity, i.e. the cycle
     \[ \Delta_X - \Xi\circ \Phi\ \ \ \in\ A^{\dim X}(X\times X) \]
     is homologically trivial.  
      
      We now consider things family-wise, i.e. we construct universal families $\XX\to B$ and $\Ss\to B$, where 
       \[ B\ \ \subset\ \PP H^0\bigl(Y_1\times Y_2,\OO_{Y_1\times Y_2}(1,1)\bigr) \]
       parametrizes all divisors $H$ such that both $X:=X_H$ and $S:=S_H$ are smooth and dimensionally transverse.
     
    Applying Voisin's Hilbert schemes argument  \cite[Proposition 3.7]{V0} (cf. also \cite[Proposition 2.11]{Lacub}) to this set-up, we may assume that the correspondences $\Phi$ and $\Xi$ are generically defined (with respect to $B$), and so in particular
      \[      \Delta_X - \Xi\circ \Phi\ \ \ \in\ GDA^{\dim X}(X\times X) \ .\]
  We observe that $H^{\dim X}_{tr}(X,\QQ)\cong H^{\dim S}_{tr}(S,\QQ)$ (this follows from the decomposition \eqref{so}), and so $H^{\dim X}_{tr}(X,\QQ)\not=0$; all conditions of Proposition \ref{spread} are fulfilled.
  Applying Proposition \ref{spread} to the cycle $  \Delta_X - \Xi\circ \Phi$, we find that a modification of this cycle vanishes:
   \[      \Delta_X - \Xi\circ \Phi -\gamma=0\ \ \ \hbox{in}\  A^{\dim X}(X\times X) \ ,\]    
   where 
   \[ \gamma\in \Bigl\langle (p_1)^\ast GDA^\ast(X), (p_2)^\ast GDA^\ast(X)\Bigr\rangle\] 
   is a decomposable cycle.
   This translates into the fact that (up to adding some trivial motives $\one(\ast)$ and modifying the correspondences $\Phi$ and $\Xi$) the composition
   \[  h(X) \   \xrightarrow{\Phi}\ h^{}_{}(S)(-m) \oplus \bigoplus \one(\ast)\ \xrightarrow{\Xi}\  h(X)\ \ \ \hbox{in}\ \MM_{\rm rat} \]
   is the identity, which proves the proposition.

    (Finally, the statement in parentheses is a straightforward consequence of the injection of motives: taking Chow groups, one obtains an injection
      \[ A^j_{hom}(X)\ \hookrightarrow\ A^{j-m}_{hom}(S)\ .\]
      But the group on the right vanishes for $j-m>\dim S$, which means $j>  {1\over 2}(\dim X+\dim S)$.) 
  \end{proof}
   
   \begin{example} Here is a sample application of Theorem \ref{hpd}. Let $Y_1=Y_2=\Gr(2,5)\subset\PP^9$. Then $Y_2^\vee=\Gr(2,5)\subset(\PP^9)^\vee$ and $S:=Y_1\cap (f_H)^{-1}(Y_2^\vee)$ is 3-dimensional (for $H$ sufficiently general). We consider the 11-dimensional variety
   \[ X:= \bigl(\Gr(2,5)\times \Gr(2,5)\bigr)\cap H\ \ \ \subset \PP^9\times \PP^9 \ ,\]
 where $H$ is a general $(1,1)$-divisor. This $X$ is a Fano variety of Calabi--Yau type, considered in \cite[Section 3.3]{IM0}.
   Theorem \ref{hpd} implies that one has
     \[ A^j_{hom}(X)=0\ \ \ \forall\ j> 7\ ,\]
     i.e. $X$ has $\hbox{Niveau}(A^\ast(X))\le 3$ in the sense of \cite{moi}.
       \end{example}

  \section{Main result}
  
 This section contains the proof of our main result, which is as follows:
  
  \begin{theorem}\label{main} Let $X\subset U$ be the inclusion of a Fano variety $X$ in its ambient space $U$, where $X,U$ are as in Table \ref{table:1}. Then $X$ has an MCK decomposition. The Chern classes $c_j(X)$, and the image $\ima\bigl( A^\ast(U)\to A^\ast(X)\bigr)$, lie in $A^\ast_{(0)}(X)$. 
  %The correspondence $\Gamma$ linking $X$ and the associated K3 surface lies in $A^\ast_{(0)}(X\times S)$.
\end{theorem}

\subsection{A criterion}
To prove Theorem \ref{main}, we will use the following general criterion:

\begin{proposition}\label{crit} Let $\XX\to B$ be a family of smooth projective varieties. Assume the following conditions:

\noindent
(c0) each fiber $X$ has dimension $2d \ge 8$;

\noindent
(c1) each fiber $X$ has a self-dual CK decomposition $\{\pi^\ast_X\}$ which is generically defined (with respect to $B$), and $h^j(X)\cong\oplus \one(\ast)$ for $j\not=2d$;

\noindent
(c2) there exists a family of surfaces $\Ss\to B^\circ$ where $B^\circ\subset B$ is a countable intersection of non-empty Zariski opens, and for each 
$b\in B^\circ$ there is a split injection of motives
  \[ h(X_b)  \ \hookrightarrow\ h(S_b)(1-d)\oplus \bigoplus \one(\ast)\ \ \ \hbox{in}\ \MM_{\rm rat}\ .\]
  
\noindent
(c3) the family $\Ss\times_{B^\circ} \Ss\to B^\circ$ has the Franchetta property.

Then for each fiber $X$, $\{\pi^\ast_X\}$ is an MCK decomposition, and $GDA^\ast(X)\subset A^\ast_{(0)}(X)$.

\medskip\noindent
Moreover, condition (c3) may be replaced by the following:

\noindent
(c3$^\prime$) $\Ss\to B^\circ$ is a family of K3 surfaces, which is the universal family of smooth sections of a direct sum of very ample line bundles on some smooth projective ambient space $V$ with trivial Chow groups (i.e. $A^\ast_{hom}(V)=0$), and $\Ss\to B^\circ$ has the Franchetta property.
\end{proposition}

\begin{proof} Using Voisin's Hilbert schemes argument  \cite[Proposition 3.7]{V0} (cf. also \cite[Proposition 2.11]{Lacub}), one may assume  that the split injection of (c2) is generically defined (with respect to $B^\circ$). This means that there exists a relative correspondence
       $\Phi$ fitting into a commutative diagram
       \[  \begin{array}[c]{ccc}     A^j(\XX) & \xrightarrow{  \Phi_\ast} & A^{j+1-d}(\Ss)\oplus \bigoplus A^\ast(B^\circ)\\
       &&\\
       \downarrow&&\downarrow\\
       &&\\
        A^j(X) & \xrightarrow{ (\Phi\vert_b)_\ast } & \ A^{j+1-d}(S)\oplus \QQ^s  ,\\ 
        \end{array}\]
        where vertical arrows are restrictions to a fiber, and the lower horizontal arrow is induced by the injection of (c2). The same then applies to $X\times X$, i.e.
        there is a commutative diagram
         \[ \begin{array}[c]{ccc}     A^j(\XX\times_{B^\circ} \XX) & \xrightarrow{} & A^{j+2-2d}(\Ss\times_{B^\circ} \Ss)\oplus \bigoplus A^{\ast}(\Ss) \oplus \bigoplus A^\ast(B^\circ)  \\
       &&\\
       \downarrow&&\downarrow\\
       &&\\
        A^j(X\times X) & \hookrightarrow& \ A^{j+2-2d}(S\times S)\oplus \bigoplus A^{\ast}(S) \oplus \bigoplus \QQ^s,\\ 
        \end{array}\]
        where the lower horizontal arrow is split injective thanks to (c2). That is, there is an injection
        \[ GDA^j_{B^\circ}(X\times X)\ \hookrightarrow\ GDA^{j+2-2d}_{B^\circ}(S\times S)\oplus \bigoplus GDA^\ast_{B^\circ}(S)\oplus \QQ^s\ .\]
        It then follows from (c3) that $\XX\times_{B^\circ} \XX\to B^\circ$ has the Franchetta property.\footnote{(NB: in practice, one can often avoid recourse to the Hilbert scheme argument in this step. For instance, in the setting of Proposition \ref{p1} below, the split injection of (c2) is generically defined by construction, and one can apply Proposition \ref{ji2} to conclude that 
$\XX\times_{B^\circ} \XX\to B^\circ$ has the Franchetta property.)}
        
        Let us now ascertain that the CK decomposition $\{\pi^\ast_X\}$ is multiplicative. What we need to check is that for each $X=X_b$ one has
   \begin{equation}\label{this} \pi_X^k\circ \Delta_X^{sm}\circ (\pi_X^i\times \pi_X^j)=0\ \ \ \hbox{in}\ A^{4d}(X\times X\times X)\ \ \ \hbox{for\ all\ }i+j\not=k\ .\end{equation}
   A standard spread lemma (cf. \cite[Lemma 3.2]{Vo}) shows that it suffices to prove this for all $b\in B^\circ$, so we will henceforth assume that $X=X_b$ with $b\in B^\circ$.
   We note that the cycle in \eqref{this} is generically defined, and homologically trivial.
   
   Let us assume that among the three integers $(i,j,k)$, at least one is different from $2d$. Using the hypothesis $h^j(X)=\oplus\one(\ast)$ for $j\not=2d$, we find there is a (generically defined) split injection
   \[   ( \pi^{4d-i}_X\times \pi^{4d-j}_X\times\pi^k_X)_\ast A^{4d}(X\times X\times X)\ \hookrightarrow\ A^\ast(X\times X)\ .\] 
   Since
     \[  \pi_X^k\circ \Delta_X^{sm}\circ (\pi_X^i\times \pi_X^j) =  ({}^t \pi^i_X\times{}^t \pi^j_X\times\pi^k_X)_\ast (\Delta^{sm}_X) =  ( \pi^{4d-i}_X\times \pi^{4d-j}_X\times\pi^k_X)_\ast (\Delta^{sm}_X)  \]
     (where the first equality is an instance of Lieberman's lemma), the required vanishing \eqref{this} now follows from the Franchetta property for $\XX\times_B \XX\to B$.
     
     It remains to treat the case $i=j=k=2d$. Using the split injection of motives (2) and taking the tensor product, we find there is a split injection of Chow groups
     \[ A^j(X\times X\times X)\ \hookrightarrow\ A^{j+3-3d}(S^3)\oplus \bigoplus A^\ast(S^2)\oplus \bigoplus A^\ast(S)\oplus \QQ^s\ .\]
     Moreover (just as we have seen above for $X^2$), this injection respects generically defined cycles, i.e. there is an injection
     \[  GDA^j(X\times X\times X)\ \hookrightarrow\ GDA^{j+3-3d}(S^3)\oplus \bigoplus GDA^\ast(S^2)\oplus \bigoplus GDA^\ast(S)\oplus \QQ^s\ .\]     
     In particular, taking $j=4d$ we find an injection
     \[ GDA^{4d}(X\times X\times X)\ \hookrightarrow\ GDA^{d+3}(S^3)\oplus \bigoplus GDA^\ast(S^2)\oplus \bigoplus GDA^\ast(S)\oplus \QQ^s\ .\]    
     By assumption, $d\ge 4$ and so the summand $GDA^{d+3}(S^3)$ vanishes for dimension reasons.
     The required vanishing \eqref{this} then follows from the Franchetta property for $\Ss\times_{B^\circ} \Ss$. This proves that $\{\pi^\ast_X\}$ is MCK.
     
     To see that $GDA^\ast(X)\subset A^\ast_{(0)}(X)$, it suffices to note that
     \[ (\pi^k_X)_\ast GDA^j(X)\ \ \ \ (k\not=2j) \]
     is generically defined, and homologically trivial. The Franchetta property for $\XX\to B$ (which is implied by the Franchetta property for $\XX\times_{B^\circ} \XX$) then implies the vanishing
             \[ (\pi^k_X)_\ast GDA^j(X)=0\ \ \ \ (k\not=2j)\ , \]
             and so $GDA^j(X)\subset (\pi^{2j}_X)_\ast A^j(X)=:  A^j_{(0)}(X)$.
             
      Let us now proceed to show that condition (c3$^\prime$) implies condition (c3). The hypotheses of (c3$^\prime$) imply that $B^\circ$ is a Zariski open in some $\bar{B}:=\PP H^0(V,\oplus_{j=1}^s L_j)$ which is isomorphic to $\PP^r$.
            The very ampleness assumption implies that
      \[ \pi\colon\ \  \Ss\times_{\bar{B}} \Ss\ \to\ V\times V \]
      is a $\PP^{r-2s}$-bundle over $(V\times V)\setminus \Delta_V$ and a $\PP^{r-s}$-bundle over $\Delta_V$. That is, $\pi$ is a {\em stratified projective bundle\/} in the sense of \cite{FLV}. As such, \cite[Proposition 5.2]{FLV} implies the equality
      \[  GDA^\ast(S\times S)= \ima\Bigl( A^\ast(V\times V)\to A^\ast(S\times S)\Bigr) +  \Delta_\ast GDA^\ast(S)\ ,\]
      where $\Delta\colon S\to S\times S$ is the inclusion along the diagonal. Since $V$ has trivial Chow groups, one has $A^\ast(V\times V)\cong A^\ast(V)\otimes A^\ast(V)$. Moreover,
      $\Ss\to V$ is a projective bundle and so \cite[Proposition 5.2]{FLV} gives $GDA^\ast(S)=\ima\bigl (A^\ast(V)\to A^\ast(S)\bigr)$. It follows that the above equality reduces to
      \begin{equation}\label{gda}  GDA^\ast(S\times S)=\Bigl\langle (p_1)^\ast GDA^\ast(S), (p_2)^\ast GDA^\ast(S), \Delta_S\Bigr\rangle\ \end{equation}    
      (where $p_1, p_2$ denote the projection from $S\times S$ to first resp. second factor). By assumption, $S$ is a K3 surface and the Franchetta property holds for $\Ss\to B^\circ$, which means that
      \[  GDA^\ast(S) = \QQ \oplus GDA^1(S) \oplus \QQ[o]\ ,\]
      where $o\in A^2(S)$ is the Beauville--Voisin class \cite{BV}.
      Given a divisor $D\in A^1(S)$, it is known that
      \[  \Delta_S\cdot (p_j)^\ast(D)=\Delta_\ast(D)= D\times o + o\times D\ \ \ \hbox{in}\ A^3(S\times S)\ \]
      \cite[Proposition 2.6(a)]{BV}. Also, it is known that
      \[ \Delta_S\cdot (p_j)^\ast(o)= \Delta_\ast(o)= o\times o\ \ \ \hbox{in}\ A^4(S\times S) \]
        \cite[Proposition 2.6(b)]{BV}. It follows that the right-hand side of \eqref{gda} is {\em decomposable\/} in codimension $>2$, i.e.
        \[ \begin{split}  \Bigl\langle (p_1)^\ast GDA^\ast(S), (p_2)^\ast GDA^\ast(S), \Delta_S\Bigr\rangle  \cap A^j(S\times S) &=\\\Bigl\langle (p_1)^\ast GDA^\ast(S), (p_2)^\ast GDA^\ast(S)\Bigr\rangle\cap A^j(S\times S)& \ \ \ \ \ \ \forall j\not=2\ .\\ \end{split}\]
        Since we know that $GDA^\ast(S)$ injects into cohomology (this is the Franchetta property for $\Ss\to B^\circ$), equality \eqref{gda}
        (plus the K\"unneth decomposition in cohomology) now implies that 
          \[ GDA^j(S\times S)\ \to\ H^{2j}(S\times S,\QQ) \] is injective for $j\not=2$. 
          
          For the case $j=2$, it suffices to remark that $\Delta_S$ is linearly independent from the decomposable part in cohomology (for otherwise $H^{2,0}(S)$ would be zero, which is absurd). The injectivity of 
          \[ GDA^2(S\times S)\ \to\  H^4(S\times S,\QQ) \] 
          then follows from \eqref{gda} plus the injectivity of $GDA^\ast(S)\to H^\ast(S,\QQ)$. This shows that condition (c3$^\prime$) implies condition (c3); the proposition is proven.
          \end{proof}

\subsection{Verifying the criterion: part 1}

\begin{proposition}\label{p1} The following families verify the conditions of Proposition \ref{crit}: the universal families $\XX\to B$ of Fano varieties of type M1, M6, M7, M8, M9, M10.
\end{proposition}

\begin{proof} The existence of a generically defined CK decomposition is an easy consequence of the fact that the Fano varieties $X$ under consideration are complete intersections in an ambient space $U$ with trivial Chow groups, cf. for instance \cite[Lemma 3.6]{V0}. This takes care of conditions (c0) and (c1) of Proposition \ref{crit}.

To verify condition (c2), we use Cayley's trick (Theorem \ref{ji}). The K3 surface $S$ associated to the Fano variety $X$ is a complete intersection in an ambient space $V$ as indicated in Table \ref{table:2}. Let us write $2d:=\dim X$. The ambient spaces $V$ that occur all have trivial Chow groups, and so Theorem \ref{ji}
gives the split injection of motives
  \[ h(X)\ \hookrightarrow\ h(S)(1-d)\oplus \bigoplus\one(\ast)\ \ \ \hbox{in}\ \MM_{\rm rat}\ ,\]
  i.e. condition (c2) is verified.

\medskip  
\begin{table}[h]
\centering
\begin{tabular}{||c c c  c||} 
 \hline
 ${\hbox{Label\ in\ \cite{FM}}}$ & $X$ & $\dim X$ &     $ S\subset V$   \\ 
 [0.5ex] 
 \hline\hline
  M1 & $X_{(1,1,1)}\subset\PP^3\times\PP^3\times\PP^3$ & 8 &  $S_{1^4}\subset \PP^3\times\PP^3$    \\ 
 M6 & $X_{(1,1)}\subset \SSS_5\times\PP^7$ &  16 & $S_{1^8}\subset \SSS_5$  \\
 M7 &   $X_{(1,1)}\subset\Gr(2,6)\times \PP^5$ & 12 &  $S_{1^6}\subset \Gr(2,6)$  \\ 
 M8 &  $X_{(1,1)}\subset\SGr(2,6)\times \PP^4$ & 10 &  $S_{1^5}\subset \SGr(2,6)$  \\ 
 M9 & $X_{(1,1)}\subset S_2 \Gr(2,6)\times \PP^3$ & 8 &  $S_{1^4}\subset S_2 \Gr(2,6)$   \\  
 M10 & $X_{(1,1)}\subset\SGr(3,6)\times \PP^3$ & 8 &  $S_{1^4}\subset    \SGr(3,6)$  \\  
  [1ex] 
 \hline
\end{tabular}
\caption{Fano varieties $X$ and their associated K3 surface $S$.} 
\label{table:2}
\end{table}
\medskip

We observe that all ambient spaces $V$ in Table \ref{table:2} have trivial Chow groups. To verify condition (c3$^\prime$), it only remains to check the Franchetta property for the families $\Ss\to B^\circ$. In all these cases, $\Ss\to V$ is a projective bundle, and so (using the projective bundle formula, or lazily applying \cite[Proposition 5.2]{FLV}) we find equality
  \[ GDA^j_{B^\circ}(S)=\ima\bigl( A^j(V)\to A^j(S)\bigr)\ .\]
  Let us check that the right-hand side injects into cohomology. This is non-trivial only in codimension $j=2$. For the family M1, it suffices to observe that
  $A^2(\PP^3\times\PP^3)$ is generated by intersections of divisors, and so 
    \[ \ima\Bigl( A^2(\PP^3\times\PP^3)\to A^2(S)\Bigr) = \QQ[o] \]
    injects into cohomology. For the family M7, it suffices to check that the restriction of $c_2(Q)\in A^2(\Gr(2,6))$ (where $Q$ denotes the universal quotient bundle) to $S$ is proportional to $o$; this is done in \cite[Proposition 2.1]{PSY}. For the family M6,  we may as well verify that
    \[ \ima \Bigl(A^2(\OGr(5,10))\to A^2(S)\Bigr)=\QQ[o] \]
    (recall that $\SSS_5$ is a connected component of $\OGr(5,10)$ in its spinor embedding), this is taken care of in \cite[Proposition 2.1]{PSY}.
For the families M8 and M9, since $\SGr(2,6)$ and $S_2 \Gr(2,6)$ are complete intersections (of dimension 7 resp. 6) inside $\Gr(2,6)$, there is an isomorphism 
   \[ A^2(S_2 \Gr(2,6))\ \xrightarrow{\cong}\  A^2(\SGr(2,6))\  \xrightarrow{\cong}\   A^2(\Gr(2,6))\ .\]
   The case M7 then guarantees that $\ima\bigl( A^2 (V)\to A^2(S)\bigr)$ is spanned by $o$. Finally, for the case M10 one observes that $A^2(\SGr(3,6))\cong\QQ$ (this follows from \cite[Proposition 2.1]{vdG}, where $\SGr(3,6)$ is denoted $Y_3$), and so
   \[ \ima  \Bigl(A^2(\SGr(3,6))\to A^2(S)\Bigr)=\QQ[h^2]= \QQ[o] \ .\]  
   \end{proof}
 
 \begin{remark} It seems likely that the families M7, M8, M9, M10 can be related to one another via (a higher-codimension version of) the game of {\em projections\/} and {\em jumps\/} of \cite[Sections 3.3 and 3.4]{BFM}. This might simplify the above argument.
 \end{remark}

 \subsection{Verifying the criterion: part 2}
 
 \begin{proposition}\label{p2} The following families verify the conditions of Proposition \ref{crit}: the universal families $\XX\to B$ of Fano varieties of type 
 M3 and M4.
 \end{proposition}
 
 \begin{proof} The existence of a generically defined CK decomposition follows as above. The difference with the above is that the families M3 and M4 are {\em not\/} in the form of Cayley's trick; hence, to check condition (c2) we now apply Theorem \ref{hpd} rather than Theorem \ref{ji}.
 
 For the case M3, Theorem \ref{hpd} applies with $Y_1=\Gr(2,5)$ and $Y_2=Q_5$ a 5-dimensional quadric embedded in $\PP^9$. Let $B^\circ\subset B$ be the open parametrizing Fano varieties $X$ of type M3 for which, in the notation of Theorem \ref{hpd}, $S$ is a smooth surface.
 For each $X=X_b$ with $b\in B^\circ$, Theorem \ref{hpd} gives an injection of motives
   \[ h(X)\ \hookrightarrow\ h(S)(-4) \oplus \bigoplus\one(\ast)\ \ \ \hbox{in}\ \MM_{\rm rat}\ .\]
 Since quadrics are projectively self-dual, this $S$ is the intersection of $\Gr(2,5)$ with a quadric and 3 hyperplanes in $\PP^9$; this is Mukai's model for the general K3 surface of genus 6.
 That the family $\Ss\to B^\circ$ has the Franchetta property is proven in \cite{PSY}. This takes care of conditions (c2) and (c3) of Proposition \ref{crit}.
 
 For the family M4, Theorem \ref{hpd} applies again, with $Y_1=\SGr(2,5)$ and $Y_2=Q_4$ a 4-dimensional quadric embedded in $\PP^9$.
 Note that $Y_1$ is a hyperplane section of $\Gr(2,5)$ under its Pl\"ucker embedding. Again, let $B^\circ\subset B$ denote the open where both $X$ and $S$ are smooth dimensionally transverse. Theorem \ref{hpd} now gives an injection of motives
   \[ h(X)\ \hookrightarrow\ h(S)(-3) \oplus \bigoplus\one(\ast)\ \ \ \hbox{in}\ \MM_{\rm rat}\ ,\] 
   where $S$ is again the intersection of $\Gr(2,5)$ with a quadric and 3 hyperplanes in $\PP^9$. The family $\Ss\to B^\circ$ is now the family of all smooth 2-dimensional complete intersections of $\SGr(2,5)$ with a quadric and 2 hyperplanes. One has that $\Ss\to \SGr(2,5)$ is a projective bundle, and so (as before)
   \[ GDA^2(S)=\ima\Bigl( A^2(\SGr(2,5))\to A^2(S)\Bigr)\ .\]
 But $A^2(\Gr(2,5))\to A^2(\SGr(2,5))$ is an isomorphism (weak Lefschetz), and so $GDA^2(S)=\QQ[o]$ as for the family M3. All conditions of Proposition \ref{crit} are verified. 
    \end{proof}

 \subsection{Proof of theorem}
 
 \begin{proof}(of Theorem \ref{main}) For the families B1 and B2 the result was proven in \cite{40}. The family S2 was treated in \cite{S2}. For the remaining families, we have checked (Propositions \ref{p1} and \ref{p2}) that Proposition \ref{crit} applies, which gives a generically defined MCK decomposition. The Chern classes $c_j(X)$, as well as the image
 $\ima\bigl(A^\ast(U)\to A^\ast(X)\bigr)$, are clearly generically defined, and so they are in $A^\ast_{(0)}(X)$ thanks to Proposition \ref{crit}.
  \end{proof}

 \section{A consequence}

 \begin{corollary}\label{cor} Let $X\subset U$ be the inclusion of a Fano variety $X$ in its ambient space $U$, where $X,U$ are as in Table \ref{table:1}. Let $\dim X=2d$. Let $R^\ast(X)\subset A^\ast(X)$ be the $\QQ$-subalgebra
  \[ R^\ast(X):=\Bigl\langle A^1(X), A^2(X), \ldots, A^d(X), c_j(X),\ima\bigl(A^\ast(U)\to A^\ast(X)\bigr)\Bigr\rangle\ \ \ \subset A^\ast(X)\ .\]
  Then $R^\ast(X)$ injects into cohomology under the cycle class map. 
 \end{corollary}
 
 \begin{proof} This is a formal consequence of the MCK paradigm. We know (Theorem \ref{main}) that $X$ has an MCK decomposition, and $c_j(X)$ and $\ima\bigl(A^\ast(U)\to A^\ast(X)\bigr)$ are in $A^\ast_{(0)}(X)$. Moreover, we know that
   \begin{equation}\label{vani} A^j_{hom}(X) =0\ \ \ \ \forall j\not=d+1 \end{equation}
   (indeed, the injection of motives of Proposition \ref{crit}(c2) induces an injection $A^j_{hom}(X)\hookrightarrow A^{j+1-d}(S)$ where $S$ is a K3 surface). This means that
   \[ A^j(X) =A^j_{(0)}(X)\ \ \ \ \forall j\not= d+1\ ,\]
   and so 
   \[ R^\ast(X)\ \ \subset\ A^\ast_{(0)}(X) \ .\]
 It only remains to check that $A^\ast_{(0)}(X)$ injects into cohomology under the cycle class map. In view of \eqref{vani}, this reduces to checking that the cycle class map induces an injection
   \[ A^{d+1}_{(0)}(X)\ \ \hookrightarrow\ H^{2d+2}(X,\QQ)\ .\]
   By construction, the correspondence $\pi_X^{2d+2}$ is supported on a subvariety $V\times W\subset X\times X$, where $V,W\subset X$ are (possibly reducible) subvarieties of dimension $\dim V=d+1$ and $\dim W=d-1$. As in \cite{BS}, the action of $\pi^{2d+2}_X$ on $A^{d+1}(X)$ factors over $A^0(\wt{W})$, where $\wt{W}\to W$ is a resolution of singularities. In particular, the action of $\pi^{2d+2}_X$ on $A^{d+1}_{hom}(X)$  factors over  $A^0_{hom}(\wt{W})=0$ and so is zero. But the action of  $\pi^{2d+2}_X$ on $A^{d+1}_{(0)}(X)$  is the identity, and so
    \[     A^{d+1}_{(0)}(X)\cap A^{d+1}_{hom}(X)=0\ ,\]
    as requested.
  \end{proof}

 \vskip1cm
\begin{nonumberingt} Thanks to Lie Fu and Charles Vial for lots of enriching exchanges around the topics of this paper. Thanks to the referee for pertinent comments. Thanks to Kai who is a great expert on Harry Potter trivia.
\end{nonumberingt}

\end{document}